\documentclass[english]{article}
\usepackage[T1]{fontenc}
\usepackage[latin9]{inputenc}
\usepackage{geometry}
\usepackage{units}
\usepackage{amsmath}
\usepackage{amssymb}

\makeatletter
\usepackage{tikz}

\makeatother

\usepackage{babel}
\begin{document}

\title{On graceful labelings of trees}

\author{Edinah K. Gnang \thanks{Department of Applied Mathematics and Statistics, Johns Hopkins University,
egnang1@jhu.edu}}
\maketitle
\begin{abstract}
We prove via a composition lemma, the Kötzig-Ringel-Rosa conjecture,
better known as the Graceful Labeling Conjecture. We also prove via
a stronger version of the composition lemma a stronger form of the
Graceful Labeling Conjecture.
\end{abstract}

\section{Introduction}

The Kötzig-Ringel-Rosa \cite{R64} conjecture, better known as the
\emph{Graceful Labeling Conjecture }(GLC), asserts that every tree
admits a graceful labeling. For a detailed survey of the extensive
literature on this problem, see \cite{Gal05}. For our purpose, we
define the graceful labeling of a graph to be a vertex labeling which
results in a bijection between vertex labels and \emph{induced subtractive
edge labels}. For notational convenience, let
\[
\mathbb{Z}_{n}:=\left[0,n\right)\cap\mathbb{Z}
\]
 Induced subtractive edge labels correspond to absolute differences
of integers assigned to the vertices spanning each edge. Our discussion
is based upon a functional reformulation of the GLC. A rooted tree
on $n>0$ vertices is associated with a function
\begin{equation}
f:\mathbb{Z}_{n}\rightarrow\mathbb{Z}_{n}\:\text{ subject to }\:\left|f^{(n-1)}\left(\mathbb{Z}_{n}\right)\right|=1,\label{functional_reformulation}
\end{equation}
\[
\text{where}
\]
\[
\forall\,i\in\mathbb{Z}_{n},\;f^{(0)}\left(i\right)\,:=i,\mbox{ and }\forall\,k\ge0,\;f^{(k+1)}\left(i\right)=f^{(k)}\left(f\left(i\right)\right)=f\left(f^{(k)}\left(i\right)\right).
\]
Every $f\in\mathbb{Z}_{n}^{\mathbb{Z}_{n}}$ has a corresponding \emph{functional
directed graph }denoted \emph{$G_{f}$} whose vertex and edge sets
are respectively
\[
V\left(G_{f}\right):=\mathbb{Z}_{n}\quad\text{ and }\quad E\left(G_{f}\right):=\left\{ \left(i,f\left(i\right)\right)\,:\,i\in\mathbb{Z}_{n}\right\} .
\]
When $f$ is subject to Eq. (\ref{functional_reformulation}), the
corresponding functional directed graph $G_{f}$ is a directed rooted
tree with an additional loop edge placed at the fixed point and whose
edges are oriented to point towards the fixed point or root of the
tree. Examples of induced edge labels associated with a functional
directed graph $G_{f}$ of $f\in\mathbb{Z}_{n}^{\mathbb{Z}_{n}}$
include :
\begin{itemize}
\item Induced subtractive edge labels given by $\left\{ \left|f\left(i\right)-i\right|\,:\,i\in\mathbb{Z}_{n}\right\} $
and $G_{f}$ is graceful if there exist $\sigma\in\text{S}_{n}\subset\mathbb{Z}_{n}^{\mathbb{Z}_{n}}$
such that $\left\{ \left|\sigma f\left(i\right)-\sigma\left(i\right)\right|\,:\,i\in\mathbb{Z}_{n}\right\} =\mathbb{Z}_{n}$
.
\item More general $\tau$-induced edge labels given by $\left\{ \tau\left(i,f\left(i\right)\right)\,:\,i\in\mathbb{Z}_{n}\right\} $,
for some given function 
\[
\tau\in\mathbb{Z}_{n}^{\mathbb{Z}_{n}\times\mathbb{Z}_{n}}
\]
 and $G_{f}$ is $\tau$-Zen for some $\tau\in\mathbb{Z}_{n}^{\mathbb{Z}_{n}\times\mathbb{Z}_{n}}$
if there exists $\sigma\in\text{S}_{n}\subset\mathbb{Z}_{n}^{\mathbb{Z}_{n}}$
such that
\[
\left\{ \tau\left(\sigma f\left(i\right),\sigma\left(i\right)\right)\,:i\in\mathbb{Z}_{n}\right\} =\mathbb{Z}_{n}
\]
\end{itemize}
The following proposition expresses a necessary and sufficient condition
for a functional directed graph to be graceful.\\
\textbf{}\\
\textbf{Proposition} 0: ( \emph{Graceful expansion} ) Let $\nicefrac{\mbox{S}_{n}}{\mathcal{I}_{n}}$
denote the partition of the symmetric group into $\frac{n!}{2}$ equivalence
classes, each made up of an unordered permutation pair $\left\{ \gamma,\,\left(n-1\right)-\gamma\right\} $,
for all $\gamma\in$ S$_{n}$. Let id $\in$ S$_{n}$ denote its identity
element and $\varphi$ the involution $\left(n-1\right)-$ id. For
$f\in\mathbb{Z}_{n}^{\mathbb{Z}_{n}}$ , let Aut$G_{f}\subseteq$
$S_{n}$ denote the automorphism group of $G_{f}$. Then $G_{f}$
is graceful iff there exist a nonempty $\mathcal{G}\subset\nicefrac{\mbox{S}_{n}}{\mathcal{I}_{n}}$
as well as a sign function $\mathfrak{s}\in\left\{ -1,0,1\right\} ^{\mathcal{G}\times\mathbb{Z}_{n}}$
such that
\begin{equation}
f\left(i\right)=\sigma_{\gamma}^{(-1)}\varphi^{\left(t\right)}\left(\varphi^{\left(t\right)}\sigma_{\gamma}\left(i\right)+\left(-1\right)^{t}\cdot\mathfrak{s}\left(\gamma,\sigma_{\gamma}\left(i\right)\right)\cdot\gamma\sigma_{\gamma}\left(i\right)\right),\quad\forall\begin{array}{c}
i\in\mathbb{Z}_{n}\\
\gamma\in\mathcal{G}\,\text{ and }\,t\in\left\{ 0,1\right\} 
\end{array},\label{Graceful Expansion}
\end{equation}
for some coset representative $\sigma_{\gamma}\in\nicefrac{\text{S}_{n}}{\text{Aut}G_{f}}.$
The index notation for $\sigma_{\gamma}$ emphasizes the dependence
of the coset representative on the permutation parameter $\gamma$.
As an illustration consider
\[
f:\mathbb{Z}_{4}\rightarrow\mathbb{Z}_{4}\;\text{s.t.}\;f\left(i\right)=\begin{cases}
\begin{array}{cc}
0 & \text{if }i=0\\
i-1 & \text{otherwise}
\end{array}\forall\,i\in\mathbb{Z}_{4},\end{cases}
\]
\[
\varphi\left(i\right)=3-i,\;\forall\,i\in\mathbb{Z}_{4}.
\]
$\mathcal{G}$ is given by 
\[
\mathcal{G}=\left\{ \gamma,\gamma^{\prime}\right\} \;\text{s.t.}\;\begin{array}{ccc}
\gamma\left(0\right) & = & 0\\
\gamma\left(1\right) & = & 2\\
\gamma\left(2\right) & = & 1\\
\gamma\left(3\right) & = & 3
\end{array}\:\text{ and }\:\begin{array}{ccc}
\gamma^{\prime}\left(0\right) & = & 3\\
\gamma^{\prime}\left(1\right) & = & 1\\
\gamma^{\prime}\left(2\right) & = & 0\\
\gamma^{\prime}\left(3\right) & = & 2
\end{array},
\]
the corresponding sign assignments are
\[
\mathfrak{s}:\mathcal{G}\times\mathbb{Z}_{4}\rightarrow\left\{ -1,0,1\right\} \;\text{s.t.}\;\begin{array}{ccc}
\mathfrak{s}\left(\gamma,0\right) & = & 0\\
\mathfrak{s}\left(\gamma,1\right) & = & 1\\
\mathfrak{s}\left(\gamma,2\right) & = & -1\\
\mathfrak{s}\left(\gamma,3\right) & = & -1
\end{array}\:\text{ and }\:\begin{array}{ccc}
\mathfrak{s}\left(\gamma^{\prime},0\right) & = & 1\\
\mathfrak{s}\left(\gamma^{\prime},1\right) & = & 1\\
\mathfrak{s}\left(\gamma^{\prime},2\right) & = & 0\\
\mathfrak{s}\left(\gamma^{\prime},3\right) & = & -1
\end{array}.
\]
Finally the coset representative permutation choices are
\[
\begin{array}{ccc}
\sigma_{\gamma}\left(0\right) & = & 0\\
\sigma_{\gamma}\left(1\right) & = & 3\\
\sigma_{\gamma}\left(2\right) & = & 1\\
\sigma_{\gamma}\left(3\right) & = & 2
\end{array}\:\text{ and }\:\begin{array}{ccc}
\sigma_{\gamma^{\prime}}\left(0\right) & = & 2\\
\sigma_{\gamma^{\prime}}\left(1\right) & = & 1\\
\sigma_{\gamma^{\prime}}\left(2\right) & = & 3\\
\sigma_{\gamma^{\prime}}\left(3\right) & = & 0
\end{array}.
\]
From which we see that for all $i\in\mathbb{Z}_{4}$, $t\in\left\{ 0,1\right\} $
and $\gamma\in\mathcal{G}$ 
\[
f\left(i\right)=
\]
\[
\sigma_{\gamma}^{\left(-1\right)}\varphi^{\left(t\right)}\left(\varphi^{\left(t\right)}\sigma_{\gamma}\left(i\right)+\left(-1\right)^{t}\cdot\mathfrak{s}\left(\gamma,\sigma_{\gamma}\left(i\right)\right)\cdot\gamma\sigma_{\gamma}\left(i\right)\right)=\sigma_{\gamma^{\prime}}^{\left(-1\right)}\varphi^{\left(t\right)}\left(\varphi^{\left(t\right)}\sigma_{\gamma^{\prime}}\left(i\right)+\left(-1\right)^{t}\cdot\mathfrak{s}\left(\gamma^{\prime},\sigma_{\gamma}\left(i\right)\right)\cdot\gamma^{\prime}\sigma_{\gamma^{\prime}}\left(i\right)\right).
\]
\\
\\
\emph{Proof} : On the one hand, the proof of necessity, follows from
the observation that $G_{f}$ being graceful implies the existence
of some coset representative $\sigma_{\gamma}\in\nicefrac{\text{S}_{n}}{\text{Aut}G_{f}}$
such that 
\[
\left\{ \left|\varphi^{\left(t\right)}\sigma_{\gamma}f\sigma_{\gamma}^{(-1)}\left(j\right)-\varphi^{\left(t\right)}\left(j\right)\right|\,:\,j\in\mathbb{Z}_{n}\right\} =\mathbb{Z}_{n}.
\]
Consequently, there exists a non empty $\mathcal{G}\subset\nicefrac{\mbox{S}_{n}}{\mathcal{I}_{n}}$
and a unique corresponding sign function $\mathfrak{s}\in\left\{ -1,0,1\right\} ^{\mathcal{G}\times\mathbb{Z}_{n}}$
such that 
\[
\left(\varphi^{\left(t\right)}\sigma_{\gamma}f\sigma_{\gamma}^{(-1)}\left(j\right)-\varphi^{\left(t\right)}\left(j\right)\right)=\left(-1\right)^{t}\cdot\mathfrak{s}\left(\gamma,j\right)\cdot\gamma\left(j\right),\quad\forall\begin{array}{c}
j\in\mathbb{Z}_{n}\\
\gamma\in\mathcal{G}\,\text{ and }\,t\in\left\{ 0,1\right\} 
\end{array}.
\]
\[
\implies f\sigma_{\gamma}^{(-1)}\left(j\right)=\sigma_{\gamma}^{(-1)}\varphi^{\left(-t\right)}\left(\varphi^{\left(t\right)}\left(j\right)+\left(-1\right)^{t}\cdot\mathfrak{s}\left(\gamma,j\right)\cdot\gamma\left(j\right)\right),\quad\forall\begin{array}{c}
j\in\mathbb{Z}_{n}\\
\gamma\in\mathcal{G}\,\text{ and }\,t\in\left\{ 0,1\right\} 
\end{array}.
\]
The change of variable $j=\sigma_{\gamma}\left(i\right)$ yields
\[
f\left(i\right)=\sigma_{\gamma}^{(-1)}\varphi^{\left(-t\right)}\left(\varphi^{\left(t\right)}\sigma_{\gamma}\left(i\right)+\left(-1\right)^{t}\cdot\mathfrak{s}\left(\gamma,\sigma_{\gamma}\left(i\right)\right)\cdot\gamma\sigma_{\gamma}\left(i\right)\right),\quad\forall\begin{array}{c}
i\in\mathbb{Z}_{n}\\
\gamma\in\mathcal{G}\,\text{ and }\,t\in\left\{ 0,1\right\} 
\end{array}.
\]
The equality $\varphi^{\left(-1\right)}=\varphi$ yields 
\[
f\left(i\right)=\sigma_{\gamma}^{(-1)}\varphi^{\left(t\right)}\left(\varphi^{\left(t\right)}\sigma_{\gamma}\left(i\right)+\left(-1\right)^{t}\cdot\mathfrak{s}\left(\gamma,\sigma_{\gamma}\left(i\right)\right)\cdot\gamma\sigma_{\gamma}\left(i\right)\right),\quad\forall\begin{array}{c}
i\in\mathbb{Z}_{n}\\
\gamma\in\mathcal{G}\,\text{ and }\,t\in\left\{ 0,1\right\} 
\end{array}.
\]
as claimed for the necessity direction.\\
On the other hand the proof of sufficiency follows from the observation
that $f$ admitting an expansion of the form prescribed in Eq. (\ref{Graceful Expansion})
implies that $G_{f}$ is isomorphic to the gracefully labeled functional
directed  graphs $G_{\sigma f\sigma^{(-1)}}$ and thereby completes
the proof.\\
\\
Prop. (0) motivates the notion of \emph{graceful expansion} of a function
$f\in\mathbb{Z}_{n}^{\mathbb{Z}_{n}}$ relative to the permutation
parameter $\gamma$ expressed as
\[
f\left(i\right)=\sigma_{\gamma}^{(-1)}\varphi^{\left(t\right)}\left(\varphi^{\left(t\right)}\sigma_{\gamma}\left(i\right)+\left(-1\right)^{t}\cdot\mathfrak{s}\left(\gamma,\sigma_{\gamma}\left(i\right)\right)\cdot\gamma\sigma_{\gamma}\left(i\right)\right),\:\forall\begin{array}{c}
i\in\mathbb{Z}_{n}\\
\gamma\in\mathcal{G}\,\text{ and }\,t\in\left\{ 0,1\right\} 
\end{array}.
\]
\[
\mbox{where}
\]
\[
\varphi=\left(n-1\right)-\text{id},\ \gamma\in\mathcal{G}\subset\nicefrac{\mbox{S}_{n}}{\mathcal{I}_{n}},\:\mathfrak{s}\in\left\{ -1,\,0,\,1\right\} ^{\mathcal{G}\times\mathbb{Z}_{n}}\mbox{ and }\sigma_{\gamma}\in\nicefrac{\text{S}_{n}}{\text{Aut}G_{f}}.
\]
Let GrL$\left(G_{f}\right)$ denote the subset of distinct functional
directed graphs isomorphic to $G_{f}$ whose labeling is already graceful.
Formally we write 
\[
\text{GrL}\left(G_{f}\right)\,:=\left\{ G_{\sigma f\sigma^{\left(-1\right)}}:\sigma\in\nicefrac{\text{S}_{n}}{\text{Aut}G_{f}}\;\text{ and }\;\mathbb{Z}_{n}=\left\{ \left|\sigma f\sigma^{\left(-1\right)}-j\right|\,:\,j\in\mathbb{Z}_{n}\right\} \right\} .
\]
\emph{The induced subtractive edge label sequence} of a functional
directed graph refers to the non-decreasing sequence of induced subtractive
edge labels. For instance the function in Figure 1 
\[
f\,:\,\mathbb{Z}_{6}\rightarrow\mathbb{Z}_{6}
\]
\[
\text{defined by }
\]
\[
f\left(0\right)=0,\,f\left(1\right)=0,\,f\left(2\right)=0,\,f\left(3\right)=0,\,f\left(4\right)=3,\,f\left(5\right)=3,
\]
is a connected functional directed spanning subtree of the complete
graph ( or \emph{functional tree} for short ) on $6$ vertices. The
attractive fixed point condition from Eq. (\ref{functional_reformulation})
is met since $f^{(2)}\left(\mathbb{Z}_{6}\right)=\left\{ 0\right\} $.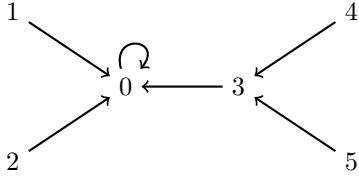
\begin{figure}
        \begin{tikzpicture}
        \node (1) at (-0.5,1) {1};
        \node (2) at (-0.5,-1) {2};
        \node (0) at (1,0) {0};
        \node (3) at (2.5,0) {3};
        \node (4) at (4,1) {4};
        \node (5) at (4,-1) {5};

        \foreach \x/\y in {1/0,2/0,3/0,4/3,5/3} {
                \draw[thick,->] (\x)--(\y);
        }
        \draw[thick,->] (0) edge [out=105,in=50,looseness=5] (0);
        \end{tikzpicture}
        \centering
        \caption{A functional directed  graph on 6 vertices.}
\end{figure} The edge set of $G_{f}$ is
\[
E\left(G_{f}\right)=\left\{ \left(0,0\right),\,\left(1,0\right),\,\left(2,0\right),\,\left(3,0\right),\,\left(4,3\right),\,\left(5,3\right)\right\} .
\]
The corresponding induced subtractive edge label sequence is
\[
\left(0,1,1,2,2,3\right).
\]
The GLC is easily verified for the families of functional star trees
associated with identically constant functions. This is seen from
the fact that functional directed graphs of identically zero functions
\[
f\,:\,\mathbb{Z}_{n}\rightarrow\mathbb{Z}_{n}
\]
\[
\text{such that}
\]
\[
f\left(i\right)=0,\quad\forall\:i\in\mathbb{Z}_{n},
\]
are gracefully labeled. For the chosen family of identically constant
functions parametrized by $n\ge2$, we have
\begin{equation}
\left|\text{GrL}\left(G_{f}\right)\right|=2.
\end{equation}
Our main results are short proofs of Composition Lemmas, from which
proofs of the GLC, the strong GLC \cite{GW18} easily follow as corollaries.

\section{The composition monoid $\mathbb{Z}_{n}^{\mathbb{Z}_{n}}$ and its
composition lemmas.}

We establish here for convenience some basic facts on the composition
monoid $\mathbb{Z}_{n}^{\mathbb{Z}_{n}}$.\\
\textbf{}\\
\textbf{Proposition} 1: ( \emph{Submonoid on functional forests} )
The subset 
\[
\left\{ h\in\mathbb{Z}_{n}^{\mathbb{Z}_{n}}:\begin{array}{c}
h\left(i\right)\le i\\
\forall\:i\in\mathbb{Z}_{n}
\end{array}\right\} 
\]
forms a submonoid of the composition monoid $\mathbb{Z}_{n}^{\mathbb{Z}_{n}}$.\\
\\
\emph{Proof} : By definition the set in question includes the identity
element. It therefore remains to show that the set is closed with
respect to composition. Let 
\[
f,g\in\left\{ h\in\mathbb{Z}_{n}^{\mathbb{Z}_{n}}:\begin{array}{c}
h\left(i\right)\le i\\
\forall\:i\in\mathbb{Z}_{n}
\end{array}\right\} 
\]
\[
\begin{array}{cc}
\implies & \forall\:i\in\mathbb{Z}_{n},\;f\left(g\left(i\right)\right)\le g\left(i\right)\le i\\
\\
\implies & fg\in\left\{ h\in\mathbb{Z}_{n}^{\mathbb{Z}_{n}}:\begin{array}{c}
h\left(i\right)\le i\\
\forall\:i\in\mathbb{Z}_{n}
\end{array}\right\} .
\end{array}
\]
Thereby completing the proof.\\
\\
Note that 
\[
n!=\left|\left\{ h\in\mathbb{Z}_{n}^{\mathbb{Z}_{n}}:\begin{array}{c}
h\left(i\right)\le i\\
\forall\:i\in\mathbb{Z}_{n}
\end{array}\right\} \right|,
\]
and functional directed graphs of functions in this particular submonoid
are each a spanning union of disjoint functional trees also called
a\emph{ functional forest}. As such, these functional forests necessarily
include at least one member of the isomorphism class of any functional
forest on $n$ vertices.\\
\\
\textbf{Proposition} 2: ( \emph{Semigroup on functional trees} ) The
set 
\[
\left\{ t\in\mathbb{Z}_{n}^{\mathbb{Z}_{n}}:\,t\left(0\right)=0\text{ and}\begin{array}{c}
t\left(i\right)<i\\
\forall\,i\in\mathbb{Z}_{n}\backslash\left\{ 0\right\} 
\end{array}\right\} ,
\]
forms a composition semigroup of the composition monoid $\left\{ h\in\mathbb{Z}_{n}^{\mathbb{Z}_{n}}:\begin{array}{c}
h\left(i\right)\le i\\
\forall\:i\in\mathbb{Z}_{n}
\end{array}\right\} $.\\
\\
\emph{Proof} :
\[
f,g\in\left\{ t\in\mathbb{Z}_{n}^{\mathbb{Z}_{n}}:\,t\left(0\right)=0\text{ and}\begin{array}{c}
t\left(i\right)<i\\
\forall\,i\in\mathbb{Z}_{n}\backslash\left\{ 0\right\} 
\end{array}\right\} 
\]
\[
\implies f\left(g\left(0\right)\right)=0.
\]
Furthermore, 
\[
\forall\:i\in\mathbb{Z}_{n}\backslash\left\{ 0\right\} ,\quad f\left(g\left(i\right)\right)=0<i\;\text{ if }\;g\left(i\right)=0\;\text{ otherwise }\ f\left(g\left(i\right)\right)<g\left(i\right)<i.
\]
Hence
\[
fg\in\left\{ t\in\mathbb{Z}_{n}^{\mathbb{Z}_{n}}:\,t\left(0\right)=0\text{ and}\begin{array}{c}
t\left(i\right)<i\\
\forall\,i\in\mathbb{Z}_{n}\backslash\left\{ 0\right\} 
\end{array}\right\} .
\]
Thus completing the proof.\\
\\
It follows as a corollary of Prop. (2) that 
\begin{equation}
\left(\underset{f\in S}{\bigcirc}f\right)\left(i\right)=0,\quad\forall\:i\in\mathbb{Z}_{n}.\label{Corr. of Prop. 2}
\end{equation}
for any multiset $S$ of elements from the semigroup
\[
\left\{ t\in\mathbb{Z}_{n}^{\mathbb{Z}_{n}}:\,t\left(0\right)=0\text{ and}\begin{array}{c}
t\left(i\right)<i\\
\forall\,i\in\mathbb{Z}_{n}\backslash\left\{ 0\right\} 
\end{array}\right\} ,
\]
such that the size of $S$ accounting for multiplicities is greater
or equal to $\left(n-1\right)$.\\
\\
\textbf{Proposition} 3: ( \emph{Submonoid sizes }) \\

1. There are at least $n!+1$ distinct submonoids of $\mathbb{Z}_{n}^{\mathbb{Z}_{n}}$
of size $n!$.

2. There are at least $n!\cdot2+n$ distinct composition submonoids
of $\mathbb{Z}_{n}^{\mathbb{Z}_{n}}$ of size $\left(n-1\right)!+1$,
if $n>3$.

3. The largest proper submonoid of $\mathbb{Z}_{n}^{\mathbb{Z}_{n}}$
has size greater or equal to $n^{n}-\frac{n!}{2}$.\\
\\
\emph{Proof} : \\
\\
1. There are $n!$ distinct submonoids of $\mathbb{Z}_{n}^{\mathbb{Z}_{n}}$
of size $n!$ parametrized by the choice of a permutation $\sigma\in$
S$_{n}\subset\mathbb{Z}_{n}^{\mathbb{Z}_{n}}$ and given by
\begin{equation}
\left\{ \sigma h\sigma^{\left(-1\right)}\in\mathbb{Z}_{n}^{\mathbb{Z}_{n}}:\begin{array}{c}
h\left(i\right)\le i\\
\forall\:i\in\mathbb{Z}_{n}
\end{array}\right\} .\label{Monoides}
\end{equation}
Note that while two submonoids parametrized by distinct permutations
$\sigma_{0}$ and $\sigma_{1}$ as described in Eq. (\ref{Monoides})
can overlap substantially, they are never equal. We see this from
the fact that the set
\[
\left\{ h\in\mathbb{Z}_{n}^{\mathbb{Z}_{n}}:\begin{array}{c}
h\left(i\right)\le i\\
\forall\:i\in\mathbb{Z}_{n}
\end{array}\right\} 
\]
includes the function 
\[
f\left(i\right)=\begin{cases}
\begin{array}{cc}
i-1 & \text{ if }i>1\\
\\
0 & \text{otherwise}
\end{array} & \forall\,i\in\mathbb{Z}_{n}\end{cases},
\]
whose functional directed graph $G_{f}$ has a trivial automorphism
group. Taken in addition with the symmetric group S$_{n}\subset\mathbb{Z}_{n}^{\mathbb{Z}_{n}}$,
these submonoids yield the desired lower bound of $n!+1$. By Cayley's
formula
\[
\left|\text{S}_{n}\cup\left(\bigcup_{\sigma\in\text{S}_{n}}\left\{ \sigma h\sigma^{\left(-1\right)}\in\mathbb{Z}_{n}^{\mathbb{Z}_{n}}:\,\begin{array}{c}
h\left(i\right)\le i\\
\forall\:i\in\mathbb{Z}_{n}
\end{array}\right\} \right)\right|=\left(n+1\right)^{\left(n-1\right)}+n!-1.
\]
2. There are $n!$ submonoids of size $\left(n-1\right)!+1$ parametrized
by the choice of a permutation $\sigma\in$ S$_{n}$ and given by
\[
\left\{ \text{id}\right\} \cup\left\{ \sigma t\sigma^{\left(-1\right)}\in\mathbb{Z}_{n}^{\mathbb{Z}_{n}}:\,t\left(0\right)=0\text{ and}\begin{array}{c}
t\left(i\right)<i\\
\forall\,i\in\mathbb{Z}_{n}\backslash\left\{ 0\right\} 
\end{array}\right\} .
\]
By invoking the argument used to prove the first part of Prop. (3),
we know that there are $\left(n-1\right)!\cdot{n \choose n-1}$ distinct
submonoids of size $\left(n-1\right)!+1$ parametrized by a choice
of an integer $j\in\mathbb{Z}_{n}$, a choice of a permutation $\sigma\in$
S$_{n}$ which fixes $j$ (i.e. $\sigma\left(j\right)=j$) and given
by
\[
\left\{ \sigma h\sigma^{\left(-1\right)}\in\mathbb{Z}_{n}^{\mathbb{Z}_{n}}:\begin{array}{c}
h\left(i\right)\le i\\
\forall\:i\in\mathbb{Z}_{n}\backslash\left\{ j\right\} 
\end{array}\text{ and }\:h\left(j\right)=j\right\} 
\]
adjoined with the constant function identically equal to $j$. Note
that, when $n>3$,
\[
n!=n\cdot\left(\left(n-1\right)!+1\right)-n\implies n!\equiv-n\mod\left(n-1\right)!+1.
\]
It follows from Lagranges' coset theorem that no submonoid of size
$\left(n-1\right)!+1$ is isomorphic to a subgroup of S$_{n}$. Furthermore,
there are $n$ special submonoids of size $\left(n-1\right)!+1$ parametrized
by the choice of an integer $j\in\mathbb{Z}_{n}$. Such a submonoid
is obtained by selecting the largest subgroup of S$_{n}$ which fixes
$j$ and adjoining to this subgroup the constant function identically
equal to $j$. The lower bound $n!\cdot2+n$ therefore follows.\\
\\
3. Consider the set
\[
\mathbb{Z}_{n}^{\mathbb{Z}_{n}}\backslash\left(\text{S}_{n}\backslash\text{A}_{n}\right),
\]
where A$_{n}$ denotes the alternating subgroup of S$_{n}$. Note
that the map
\[
h\mapsto\mathbf{A}_{G_{h}}
\]
where $h\in\mathbb{Z}_{n}^{\mathbb{Z}_{n}}$ and $\mathbf{A}_{G_{h}}\in\left\{ 0,1\right\} ^{n\times n}$
such that 
\[
\mathbf{A}_{G_{h}}\left[i,j\right]=\begin{cases}
\begin{array}{cc}
1 & \text{ if }j=h\left(i\right)\\
0 & \text{otherwise}
\end{array},\end{cases}
\]
yields a matrix representation of the composition monoid $\mathbb{Z}_{n}^{\mathbb{Z}_{n}}$.
It further follows from the rank inequality
\[
1\le\text{Rank}\left(\mathbf{A}_{G_{g}}\cdot\mathbf{A}_{G_{f}}\right)=\text{Rank}\left(\mathbf{A}_{G_{f}}\cdot\mathbf{A}_{G_{g}}\right)\le\min\left\{ \text{Rank}\left(\mathbf{A}_{G_{f}}\right),\text{Rank}\left(\mathbf{A}_{G_{g}}\right)\right\} \le n,
\]
that 
\[
\left\{ \sigma f,f\sigma\right\} \subset\mathbb{Z}_{n}^{\mathbb{Z}_{n}}\backslash\text{S}_{n},\quad\forall\;\begin{cases}
\begin{array}{ccc}
\sigma & \in & \text{S}_{n}\\
\\
f & \in & \mathbb{Z}_{n}^{\mathbb{Z}_{n}}\backslash\text{S}_{n}
\end{array}\end{cases}.
\]
Recall from Lagranges' coset theorem that no proper subgroup of S$_{n}$
is larger than the alternating group. The largest submonoid of size
less than $n^{n}$ has size at least 
\[
\left|\mathbb{Z}_{n}^{\mathbb{Z}_{n}}\backslash\left(\text{S}_{n}\backslash\text{A}_{n}\right)\right|=n^{n}-\frac{n!}{2}.
\]
Thus completing the proof.\\
\\
\textbf{}\\
\textbf{Proposition }4a: ( \emph{Preliminary Min Composition Proposition}
). For all $f\in\mathbb{Z}_{n}^{\mathbb{Z}_{n}}$, we have
\[
1\le\min_{\sigma\in\text{S}_{n}}\left|\left\{ \left|\sigma f\sigma^{\left(-1\right)}\left(i\right)-i\right|:i\in\mathbb{Z}_{n}\right\} \right|\le\rho_{f}.
\]
where $\rho_{f}$ denotes the minimum number of edge deletions required
in $G_{f}$ to obtain a spanning subgraph which is a union of disjoint
paths, free of any loop edge.\\
\\
\emph{Proof} : The lower bound is attained when the $n$ edges of
$G_{f}$ are assigned the same induced subtractive edge label. We
justify the upper-bound by considering any possible deletion of $\rho_{f}$
edges from $G_{f}$ which result in a spanning union of disjoint paths,
free of any loop edge. We then sequentially label vertices of each
path starting from one endpoint and increasing by one the label for
each vertex encountered as we move along the path towards the second
endpoint. The procedure described, greedily maximizes the number of
edges assigned the subtractive edge label equal to one. The only edges
in the proposed relabeling of $G_{f}$ whose induced subtractive edge
labels possibly different from one, correspond to the $\rho_{f}$
deleted edges. Thus completing the proof.\\
\\
Note that the lower bound is tight for any functional directed  graph
$G_{\sigma}$ associated with $\sigma\in$ S$_{n}$ which are involutions
i.e. $\sigma^{\left(-1\right)}=\sigma$ . While the upper bound is
tight for any functional directed graph which is a single spanning
directed cycle.\\
\\
\\
\textbf{Proposition} 4b: (\emph{ Preliminary Max Composition Proposition}
). For all $f\in\mathbb{Z}_{n}^{\mathbb{Z}_{n}}$, we have
\[
n-\rho_{f}+\max\left(0,\,\left(\text{\# of self-loops in }G_{f}\right)-1\right)\le\max_{\sigma\in\text{S}_{n}}\left|\left\{ \left|\sigma f\sigma^{\left(-1\right)}\left(i\right)-i\right|:i\in\mathbb{Z}_{n}\right\} \right|\le n
\]
where $\rho_{f}$ denotes the minimum number of edge deletions required
in $G_{f}$ to obtain a spanning subgraph which is a union of disjoint
paths, free of any loop edge.\\
\\
\emph{Proof} : The upper bound follows from the observation $\left|E\left(G_{f}\right)\right|=n$.
We justify the lower-bound by considering every possible deletion
of $\rho_{f}$ edges from $G_{f}$ which result in a spanning union
of disjoint paths, free of any loop edge. We then sequentially label
vertices of each path starting from one endpoint and alternating between
largest and smallest unassigned label for each vertex encountered
as we move along the path towards the second endpoint. The procedure
greedily maximizes the number of distinct edge labels. The only edges
in the proposed relabeling of $G_{f}$ whose induced subtractive edge
labels possibly repeat, correspond to the $\rho_{f}-\max\left(0,\left|\left\{ \left(j,j\right):j\in\mathbb{Z}_{n}\right\} \cap E\left(G_{f}\right)\right|-1\right)$
deleted edges. Thus completing the proof.\\
\\
Note that the upper bound is tight for any graceful functional directed
graph. The lower bound is tight for any functional directed graph
which is a single spanning directed cycle.\\
\\
\textbf{Lemma} 5: ( \emph{Composition Lemma} ) 
\[
f,g\in\left\{ t\in\mathbb{Z}_{n}^{\mathbb{Z}_{n}}:\,t\left(0\right)=0\text{ and}\begin{array}{c}
t\left(i\right)<i\\
\forall\,i\in\mathbb{Z}_{n}\backslash\left\{ 0\right\} 
\end{array}\right\} \Rightarrow\max_{\sigma\in\text{S}_{n}}\left|\left\{ \left|\sigma fg\sigma^{(-1)}\left(i\right)-i\right|:i\in\mathbb{Z}_{n}\right\} \right|\le\max_{\sigma\in\text{S}_{n}}\left|\left\{ \left|\sigma f\sigma^{(-1)}\left(i\right)-i\right|:i\in\mathbb{Z}_{n}\right\} \right|.
\]
Prior to discussing the proof we briefly introduce the notion of simultaneous
linear maps acting upon vector inputs to a given function
\[
F\,:\,\mathbb{C}^{n\times1}\times\cdots\times\mathbb{C}^{n\times1}\rightarrow\mathbb{C}
\]
 in $m$ vector variables $\left\{ \mathbf{x}_{0},\cdots,\mathbf{x}_{i},\cdots,\mathbf{x}_{m-1}\right\} $.
A simultaneous action on input vector variables of $F$ is prescribed
by effecting upon vector variables of $F$ simultaneous linear maps
\[
\mathbf{x}_{i}\mapsto\mathbf{T}_{i}\cdot\mathbf{x}_{i},\quad\forall\:0\le i<m,
\]
where each $\mathbf{T}_{i}\in\mathbb{C}^{n\times n}$ for all $0\le i<m$.
The simultaneous linear action therefore maps $F\left(\mathbf{x}_{0},\cdots,\mathbf{x}_{i},\cdots,\mathbf{x}_{m-1}\right)$
to the function $F\left(\mathbf{T}_{0}\cdot\mathbf{x}_{0},\cdots,\mathbf{T}_{i}\cdot\mathbf{x}_{i},\cdots,\mathbf{T}_{m-1}\cdot\mathbf{x}_{m-1}\right)$.\\
\emph{}\\
\emph{Proof} : We start by discussing the setting where presumably
for some element $f$, $g$ in the semigroup we have 
\[
n>\max_{\sigma\in\text{S}_{n}}\left|\left\{ \left|\sigma fg\sigma^{(-1)}\left(i\right)-i\right|:i\in\mathbb{Z}_{n}\right\} \right|.
\]
In which case consider two new functions $f^{\prime}$, $g^{\prime}$
lying in the same semigroup such that $f^{\prime}=fg$ and $g^{\prime}$
is the identically zero function. Consequently $f^{\prime}g^{\prime}$
is also identically zero and yields the inequality
\[
\max_{\sigma\in\text{S}_{n}}\left|\left\{ \left|\sigma f^{\prime}g^{\prime}\sigma^{(-1)}\left(i\right)-i\right|:i\in\mathbb{Z}_{n}\right\} \right|=n>\max_{\sigma\in\text{S}_{n}}\left|\left\{ \left|\sigma f^{\prime}\sigma^{(-1)}\left(i\right)-i\right|:i\in\mathbb{Z}_{n}\right\} \right|,
\]
which contradicts the assertion of Lem. (5) for the new functions
$f^{\prime}$, $g^{\prime}$. It is therefore necessary for Lem. (5)
to hold that
\begin{equation}
n=\max_{\sigma\in\text{S}_{n}}\left|\left\{ \left|\sigma fg\sigma^{(-1)}\left(i\right)-i\right|:i\in\mathbb{Z}_{n}\right\} \right|.\label{setting}
\end{equation}
We prove the claim of Lem. (5) by contradiction in the setting where
Eq. (\ref{setting}) holds. Assume for the sake of establishing a
contradiction that for some elements $f$, $g$ in the semigroup on
functional trees we have
\[
n=\max_{\sigma\in\text{S}_{n}}\left|\left\{ \left|\sigma fg\sigma^{(-1)}\left(i\right)-i\right|:i\in\mathbb{Z}_{n}\right\} \right|>\max_{\sigma\in\text{S}_{n}}\left|\left\{ \left|\sigma f\sigma^{(-1)}\left(i\right)-i\right|:i\in\mathbb{Z}_{n}\right\} \right|.
\]
For notational convenience, we associate with an arbitrary vector
$\mathbf{v}\in\mathbb{C}^{n\times1}$ the Vandermonde matrix
\[
\text{Vander}\left(\mathbf{v}\right)\in\mathbb{C}^{n\times n}\:\text{ s.t. }\:\text{Vander}\left(\mathbf{v}\right)\left[i,j\right]=\left(\mathbf{v}\left[i\right]\right)^{2j},\ \forall\,0\le i,j<n,
\]
The argument considers a pair of homogeneous polynomials, both of
total degree $4{n+1 \choose 2}$ in entries of $n!\cdot2$ vector
variables $\left\{ \mathbf{x}_{0},\cdots,\mathbf{x}_{n!2-1}\right\} $.
These polynomials are parametrized by $t\in\left\{ 0,1\right\} $
and given by
\[
F_{t}\left(\mathbf{x}_{0},\cdots,\mathbf{x}_{n!2-1}\right)=\sum_{\sigma\in\text{S}_{n}}\det\left(\text{Vander}\left(\mathbf{x}_{n!\cdot t+\text{lex}\sigma}\right)\right)^{2}.
\]
The map lex$:$ S$_{n}$$\rightarrow\mathbb{Z}_{n!}$ bijectively
maps S$_{n}$ to $\mathbb{Z}_{n!}$ and is prescribed by
\[
\text{lex}\left(\sigma\right)=\sum_{k\in\mathbb{Z}_{n}}\left(n-1-k\right)!\left|\left\{ \sigma\left(i\right)<\sigma\left(k\right)\,:\,i<k<n\right\} \right|,\quad\forall\ \sigma\in\text{S}_{n}.
\]
Accordingly lex$\left(\text{id}\right)=0$  and lex$\left(n-1-\text{id}\right)=n!-1$.
Consider throughout the whole discussion fixed evaluations points
\begin{equation}
\mathbf{a}_{n!\cdot t+\text{lex}\sigma}=\left(\mathbf{A}_{G_{\sigma fg^{\left(t\right)}\sigma^{\left(-1\right)}}}-\mathbf{I}_{n}\right)\left(\begin{array}{c}
r\left(c_{\sigma}+0\right)\\
\vdots\\
r\left(c_{\sigma}+i\right)\\
\vdots\\
r\left(c_{\sigma}+n-1\right)
\end{array}\right)=r\left(\begin{array}{c}
\sigma fg^{\left(t\right)}\sigma^{\left(-1\right)}\left(0\right)-0\\
\vdots\\
\sigma fg^{\left(t\right)}\sigma^{\left(-1\right)}\left(i\right)-i\\
\vdots\\
\sigma fg^{\left(t\right)}\sigma^{\left(-1\right)}\left(n-1\right)-\left(n-1\right)
\end{array}\right),\ \forall\,\begin{array}{c}
\sigma\in\text{S}_{n}\\
t\in\left\{ 0,1\right\} 
\end{array},\label{prescribed points}
\end{equation}
for symbolic variables $r$ and $\left\{ c_{\sigma}:\sigma\in\text{S}_{n}\right\} $
to be assigned real values. Throughout the discussion we evaluate
functions at points prescribed by Eq. (\ref{prescribed points}) by
assigning to the vector variable $\mathbf{x}_{n!\cdot t+\text{lex}\sigma}$
the vector $\mathbf{a}_{n!\cdot t+\text{lex}\sigma}$ for all $\sigma\in\text{S}_{n}$
and $t\in\left\{ 0,1\right\} $.
\begin{equation}
\implies F_{t}\left(\mathbf{a}_{0},\cdots,\mathbf{a}_{n!2-1}\right)=\sum_{\sigma\in\text{S}_{n}}\det\left(\text{Vander}\left(\mathbf{a}_{n!\cdot t+\text{lex}\sigma}\right)\right)^{2}=\left|\text{GrL}\left(G_{fg^{\left(t\right)}}\right)\right|\left|\text{Aut }G_{fg^{\left(t\right)}}\right|\prod_{0\le i<j<n}\left(\left(rj\right)^{2}-\left(ri\right)^{2}\right)^{2}.\label{Center Sums}
\end{equation}
The equality in Eq. (\ref{Center Sums}) stems from the observation
that each Vandermonde determinant summand expresses the product of
differences of square induced edge labels associated with pairs of
distinct edges in some relabeling of $G_{fg^{\left(t\right)}}$. On
the one hand
\[
F_{0}\left(\mathbf{a}_{0},\cdots,\mathbf{a}_{n!2-1}\right)=\sum_{\sigma\in\text{S}_{n}}\prod_{0\le i<j<n}\left(\left(\mathbf{a}_{\text{lex}\sigma}\left[j\right]\right)^{2}-\left(\mathbf{a}_{\text{lex}\sigma}\left[i\right]\right)^{2}\right)^{2},
\]
For notational convenience let
\[
\mathcal{F}_{\sigma}\left(i\right)=\left(\mathbf{a}_{\text{lex}\sigma}\left[i\right]-\mathbf{a}_{n!+\text{lex}\sigma}\left[i\right]\right),\ \forall\,\begin{array}{c}
\sigma\in\text{S}_{n}\\
i\in\mathbb{Z}_{n}
\end{array},
\]
\[
\implies F_{0}\left(\mathbf{a}_{0},\cdots,\mathbf{a}_{n!2-1}\right)=\sum_{\sigma\in\text{S}_{n}}\prod_{0\le i<j<n}\left(\left(\mathbf{a}_{n!+\text{lex}\sigma}\left[j\right]+\mathcal{F}_{\sigma}\left(j\right)\right)^{2}-\left(\mathbf{a}_{n!+\text{lex}\sigma}\left[i\right]+\mathcal{F}_{\sigma}\left(i\right)\right)^{2}\right)^{2},
\]
\[
\implies F_{0}\left(\mathbf{a}_{0},\cdots,\mathbf{a}_{n!2-1}\right)=
\]
\[
\sum_{\sigma\in\text{S}_{n}}\prod_{0\le i<j<n}\left(\left(\mathbf{a}_{n!+\text{lex}\sigma}\left[j\right]\right)^{2}-\left(\mathbf{a}_{n!+\text{lex}\sigma}\left[i\right]\right)^{2}+2\mathcal{F}_{\sigma}\left(j\right)\mathbf{a}_{n!+\text{lex}\sigma}\left[j\right]-2\mathcal{F}_{\sigma}\left(i\right)\mathbf{a}_{n!+\text{lex}\sigma}\left[i\right]+\mathcal{F}_{\sigma}\left(j\right)^{2}-\mathcal{F}_{\sigma}\left(i\right)^{2}\right)^{2}.
\]
Making use of the expansion identity
\[
\prod_{0\le i<j<n}\left(\alpha_{ij}+\beta_{ij}\right)=\left(\prod_{0\le i<j<n}\alpha_{ij}\right)+\sum_{\begin{array}{c}
k_{ij}\in\left\{ 0,1\right\} \\
0=\underset{i<j}{\prod}k_{ij}
\end{array}}\prod_{0\le i<j<n}\alpha_{ij}^{k_{ij}}\beta_{ij}^{1-k_{ij}}.
\]
we write
\[
F_{0}\left(\mathbf{a}_{0},\cdots,\mathbf{a}_{n!2-1}\right)=\sum_{\sigma\in\text{S}_{n}}\bigg[\det\text{Vander}\left(\mathbf{a}_{n!+\text{lex}\sigma}\right)+
\]
\begin{equation}
\sum_{\begin{array}{c}
k_{ij}\in\left\{ 0,1\right\} \\
0=\underset{i<j}{\prod}k_{ij}
\end{array}}\prod_{i<j}\left.\left(\left(\mathbf{a}_{n!+\text{lex}\sigma}\left[j\right]\right)^{2}-\left(\mathbf{a}_{n!+\text{lex}\sigma}\left[i\right]\right)^{2}\right)^{k_{ij}}\left(2\mathcal{F}_{\sigma}\left(j\right)\mathbf{a}_{n!+\text{lex}\sigma}\left[j\right]-2\mathcal{F}_{\sigma}\left(i\right)\mathbf{a}_{n!+\text{lex}\sigma}\left[i\right]+\mathcal{F}_{\sigma}\left(j\right)^{2}-\mathcal{F}_{\sigma}\left(i\right)^{2}\right)^{1-k_{ij}}\right]^{2}.\label{SOS}
\end{equation}
The right hand side of Eq. (\ref{SOS}) expresses a sum of square
of real numbers and as such is zero only if each summand is zero.
By our premise $F_{0}\left(\mathbf{a}_{0},\cdots,\mathbf{a}_{n!2-1}\right)=0$,
from which it follows that
\[
\forall\,\sigma\in\text{S}_{n},\:\left(-1\right)\det\text{Vander}\left(\mathbf{a}_{n!+\text{lex}\sigma}\right)=
\]
\[
\sum_{\begin{array}{c}
k_{ij}\in\left\{ 0,1\right\} \\
0=\underset{i<j}{\prod}k_{ij}
\end{array}}\prod_{i<j}\left(\left(\mathbf{a}_{n!+\text{lex}\sigma}\left[j\right]\right)^{2}-\left(\mathbf{a}_{n!+\text{lex}\sigma}\left[i\right]\right)^{2}\right)^{k_{ij}}\left(2\mathcal{F}_{\sigma}\left(j\right)\mathbf{a}_{n!+\text{lex}\sigma}\left[j\right]-2\mathcal{F}_{\sigma}\left(i\right)\mathbf{a}_{n!+\text{lex}\sigma}\left[i\right]+\mathcal{F}_{\sigma}\left(j\right)^{2}-\mathcal{F}_{\sigma}\left(i\right)^{2}\right)^{1-k_{ij}}.
\]
Squaring both sides of the equality above and summing over $\sigma\in\text{S}_{n}$
yields
\[
F_{1}\left(\mathbf{a}_{0},\cdots,\mathbf{a}_{n!2-1}\right)=
\]
\begin{equation}
\sum_{\sigma\in\text{S}_{n}}\left[\sum_{\begin{array}{c}
k_{ij}\in\left\{ 0,1\right\} \\
0=\underset{i<j}{\prod}k_{ij}
\end{array}}\prod_{i<j}\left(\left(\mathbf{a}_{n!+\text{lex}\sigma}\left[j\right]\right)^{2}-\left(\mathbf{a}_{n!+\text{lex}\sigma}\left[i\right]\right)^{2}\right)^{k_{ij}}\left(2\mathcal{F}_{\sigma}\left(j\right)\mathbf{a}_{n!+\text{lex}\sigma}\left[j\right]-2\mathcal{F}_{\sigma}\left(i\right)\mathbf{a}_{n!+\text{lex}\sigma}\left[i\right]+\mathcal{F}_{\sigma}\left(j\right)^{2}-\mathcal{F}_{\sigma}\left(i\right)^{2}\right)^{1-k_{ij}}\right]^{2}.\label{first first}
\end{equation}
Note that both sides of Eq. (\ref{first first}) express homogeneous
polynomials with integer coefficients of degree $4{n+1 \choose 2}$
evaluated at points prescribed by Eq. (\ref{prescribed points}).
Furthermore the right hand side of Eq. (\ref{first first}) expresses
the evaluation at the prescribed evaluations points of the homogeneous
polynomial
\[
G\left(\mathbf{x}_{0},\cdots,\mathbf{x}_{n!2-1}\right)=
\]
\[
\sum_{\sigma\in\text{S}_{n}}\left[\sum_{\begin{array}{c}
k_{ij}\in\left\{ 0,1\right\} \\
0=\underset{i<j}{\prod}k_{ij}
\end{array}}\prod_{i<j}\left(\left(\mathbf{x}_{n!+\text{lex}\sigma}\left[j\right]\right)^{2}-\left(\mathbf{x}_{n!+\text{lex}\sigma}\left[i\right]\right)^{2}\right)^{k_{ij}}\left(2\mathcal{G}_{\sigma}\left(j\right)\mathbf{x}_{n!+\text{lex}\sigma}\left[j\right]-2\mathcal{G}_{\sigma}\left(i\right)\mathbf{x}_{n!+\text{lex}\sigma}\left[i\right]+\mathcal{G}_{\sigma}\left(j\right)^{2}-\mathcal{G}_{\sigma}\left(i\right)^{2}\right)^{1-k_{ij}}\right]^{2},
\]
\[
\text{where }\quad\mathcal{G}_{\sigma}\left(i\right)=\left(\mathbf{x}_{\text{lex}\sigma}\left[i\right]-\mathbf{x}_{n!+\text{lex}\sigma}\left[i\right]\right),\;\forall\,\begin{array}{c}
\sigma\in\text{S}_{n}\\
k\in\mathbb{Z}_{n}
\end{array}.
\]
On the other hand recall that
\[
F_{1}\left(\mathbf{a}_{0},\cdots,\mathbf{a}_{n!2-1}\right)=\sum_{\sigma\in\text{S}_{n}}\prod_{0\le i<j<n}\left(\left(\mathbf{a}_{n!+\text{lex}\sigma}\left[j\right]\right)^{2}-\left(\mathbf{a}_{n!+\text{lex}\sigma}\left[i\right]\right)^{2}\right)^{2}.
\]
\[
\implies F_{1}\left(\mathbf{a}_{0},\cdots,\mathbf{a}_{n!2-1}\right)=\sum_{\sigma\in\text{S}_{n}}\prod_{0\le i<j<n}\left(\left(\mathbf{a}_{\text{lex}\sigma}\left[j\right]-\mathcal{F}_{\sigma}\left(j\right)\right)^{2}-\left(\mathbf{a}_{\text{lex}\sigma}\left[i\right]-\mathcal{F}_{\sigma}\left(i\right)\right)^{2}\right)^{2},
\]
\[
\implies F_{1}\left(\mathbf{a}_{0},\cdots,\mathbf{a}_{n!2-1}\right)=\sum_{\sigma\in\text{S}_{n}}\bigg[\det\text{Vander}\left(\mathbf{a}_{\text{lex}\sigma}\right)+
\]
\[
\sum_{\begin{array}{c}
k_{ij}\in\left\{ 0,1\right\} \\
0=\underset{i<j}{\prod}k_{ij}
\end{array}}\prod_{i<j}\left.\left(\left(\mathbf{a}_{\text{lex}\sigma}\left[j\right]\right)^{2}-\left(\mathbf{a}_{\text{lex}\sigma}\left[i\right]\right)^{2}\right)^{k_{ij}}\left(-2\mathcal{F}_{\sigma}\left(j\right)\mathbf{a}_{\text{lex}\sigma}\left[j\right]+2\mathcal{F}_{\sigma}\left(i\right)\mathbf{a}_{\text{lex}\sigma}\left[i\right]+\mathcal{F}_{\sigma}\left(j\right)^{2}-\mathcal{F}_{\sigma}\left(i\right)^{2}\right)^{1-k_{ij}}\right]^{2}.
\]
Note that 
\[
n>\max_{\sigma\in\text{S}_{n}}\left|\left\{ \left|\sigma f\sigma^{(-1)}\left(i\right)-i\right|:i\in\mathbb{Z}_{n}\right\} \right|\implies0=\det\text{Vander}\left(\mathbf{a}_{\text{lex}\sigma}\right),\:\forall\,\sigma\in\text{S}_{n}
\]
from which it follows that 
\[
F_{1}\left(\mathbf{a}_{0},\cdots,\mathbf{a}_{n!2-1}\right)=
\]
\begin{equation}
\sum_{\sigma\in\text{S}_{n}}\left[\sum_{\begin{array}{c}
k_{ij}\in\left\{ 0,1\right\} \\
0=\underset{i<j}{\prod}k_{ij}
\end{array}}\prod_{i<j}\left(\left(\mathbf{a}_{\text{lex}\sigma}\left[j\right]\right)^{2}-\left(\mathbf{a}_{\text{lex}\sigma}\left[i\right]\right)^{2}\right)^{k_{ij}}\left(-2\mathcal{F}_{\sigma}\left(j\right)\mathbf{a}_{\text{lex}\sigma}\left[j\right]+2\mathcal{F}_{\sigma}\left(i\right)\mathbf{a}_{\text{lex}\sigma}\left[i\right]+\mathcal{F}_{\sigma}\left(j\right)^{2}-\mathcal{F}_{\sigma}\left(i\right)^{2}\right)^{1-k_{ij}}\right]^{2}.\label{first second}
\end{equation}
Now consider the simultaneous action on $G\left(\mathbf{x}_{0},\cdots,\mathbf{x}_{n!2-1}\right)$
prescribed by maps
\begin{equation}
\mathbf{x}_{n!\cdot t+\text{lex}\sigma}\mapsto\left(\mathbf{A}_{G_{\sigma fg^{\left(1-t\right)}\sigma^{\left(-1\right)}}}-\mathbf{I}_{n}\right)\cdot\left(\mathbf{A}_{G_{\sigma fg^{\left(t\right)}\sigma^{\left(-1\right)}}}-\mathbf{I}_{n}\right)^{+}\cdot\mathbf{x}_{n!\cdot t+\text{lex}\sigma},\;\forall\:\begin{array}{c}
\sigma\in\text{S}_{n}\\
t\in\left\{ 0,1\right\} 
\end{array},\label{map}
\end{equation}
where $\left(\mathbf{A}_{G_{\sigma fg^{\left(t\right)}\sigma^{\left(-1\right)}}}-\mathbf{I}_{n}\right)^{+}$
denotes the pseudoinverse of $\left(\mathbf{A}_{G_{\sigma fg^{\left(t\right)}\sigma^{\left(-1\right)}}}-\mathbf{I}_{n}\right)$.
Evaluating the polynomial resulting from the simultaneous action on
$G\left(\mathbf{x}_{0},\cdots,\mathbf{x}_{n!2-1}\right)$ at the prescribed
points, effectively maps Eq. (\ref{first first}) to Eq. (\ref{first second}).
We point out that although signed incidence matrices  $\left(\mathbf{A}_{G_{\sigma fg^{\left(t\right)}\sigma^{\left(-1\right)}}}-\mathbf{I}_{n}\right)$
for all $\sigma\in\text{S}_{n}$ and $t\in\left\{ 0,1\right\} $ are
singular, each map in Eq. (\ref{map}) is well defined since they
share the same null space given by
\[
\text{Null Space of }\left(\mathbf{A}_{G_{\sigma fg^{\left(t\right)}\sigma^{\left(-1\right)}}}-\mathbf{I}_{n}\right)=\text{Span of }\mathbf{1}_{n\times1},\;\forall\,\begin{array}{c}
\sigma\in\text{S}_{n}\\
t\in\left\{ 0,1\right\} 
\end{array}.
\]
Furthermore, induced subtractive edge labels are invariant to vertex
label vector translations by any element in the span of $\mathbf{1}_{n\times1}$.
Now recall that the result of performing a linear transformation to
variables of a multivariate polynomial is alternatively expressed
as the result of performing some linear transformation to its coefficients.
Consequently, the simultaneous action on $G\left(\mathbf{x}_{0},\cdots,\mathbf{x}_{n!2-1}\right)$
followed by its evaluation at the prescribed points, effects in Eq.
(\ref{first first}) the involution
\begin{equation}
\mathbf{a}_{n!\cdot t+\text{lex}\sigma}\mapsto\mathbf{a}_{n!\cdot\left(1-t\right)+\text{lex}\sigma},\;\forall\:\begin{array}{c}
\sigma\in\text{S}_{n}\\
t\in\left\{ 0,1\right\} 
\end{array}.\label{Permutation}
\end{equation}
By Eq. (\ref{first first}) and Eq. (\ref{first second}), the involution
in Eq. (\ref{Permutation}) leaves the evaluation value unchanged.
Furthermore the evaluation at the prescribed point of the result of
the simultaneous action on $G\left(\mathbf{x}_{0},\cdots,\mathbf{x}_{n!2-1}\right)$
prescribed by
\[
\begin{cases}
\begin{array}{ccc}
\mathbf{x}_{\text{lex}\sigma} & \mapsto & \left(\mathbf{A}_{G_{\sigma fg\sigma^{\left(-1\right)}}}-\mathbf{I}_{n}\right)\cdot\left(\mathbf{A}_{G_{\sigma f\sigma^{\left(-1\right)}}}-\mathbf{I}_{n}\right)^{+}\cdot\mathbf{x}_{\text{lex}\sigma}\\
 & \text{and}\\
\mathbf{x}_{n!+\text{lex}\sigma} & \mapsto & \mathbf{0}_{n\times n}\cdot\mathbf{x}_{n!+\text{lex}\sigma}
\end{array} & \forall\,\sigma\in\text{S}_{n},\end{cases}
\]
also leaves the evaluation value unchanged since
\begin{equation}
G\left(\mathbf{x}_{0},\cdots,\mathbf{x}_{n!-1},\mathbf{0}_{n\times1},\cdots,\mathbf{0}_{n\times1}\right)=\sum_{\sigma\in\text{S}_{n}}\prod_{0\le i<j<n}\left(\left(\mathbf{x}_{\text{lex}\sigma}\left[j\right]\right)^{2}-\left(\mathbf{x}_{\text{lex}\sigma}\left[i\right]\right)^{2}\right)^{2}=F_{0}\left(\mathbf{x}_{0},\cdots,\mathbf{x}_{n!2-1}\right).\label{insight 1}
\end{equation}
Observing in addition that $F_{0}\left(\mathbf{a}_{0},\cdots,\mathbf{a}_{n!2-1}\right)$
equals the evaluation at the prescribed point of the result of the
simultaneous action on $G\left(\mathbf{x}_{0},\cdots,\mathbf{x}_{n!2-1}\right)$
prescribed by
\[
\begin{cases}
\begin{array}{ccc}
\mathbf{x}_{\text{lex}\sigma} & \mapsto & \mathbf{0}_{n\times n}\cdot\mathbf{x}_{\text{lex}\sigma}\\
 & \text{and}\\
\mathbf{x}_{n!+\text{lex}\sigma} & \mapsto & \left(\mathbf{A}_{G_{\sigma f\sigma^{\left(-1\right)}}}-\mathbf{I}_{n}\right)\cdot\left(\mathbf{A}_{G_{\sigma fg\sigma^{\left(-1\right)}}}-\mathbf{I}_{n}\right)^{+}\cdot\mathbf{x}_{n!+\text{lex}\sigma}
\end{array} & \forall\,\sigma\in\text{S}_{n},\end{cases}
\]
since
\[
G\left(\mathbf{0}_{n\times1},\cdots,\mathbf{0}_{n\times1},\mathbf{x}_{n!},\cdots,\mathbf{x}_{n!2-1}\right)=\sum_{\sigma\in\text{S}_{n}}\left[\sum_{\begin{array}{c}
k_{ij}\in\left\{ 0,1\right\} \\
0=\underset{i<j}{\prod}k_{ij}
\end{array}}\left(\prod_{0\le i<j<n}\left(-1\right)^{1-k_{ij}}\right)\prod_{0\le i<j<n}\left(\left(\mathbf{x}_{n!+\text{lex}\sigma}\left[j\right]\right)^{2}-\left(\mathbf{x}_{n!+\text{lex}\sigma}\left[i\right]\right)^{2}\right)\right]^{2},
\]
\begin{equation}
\implies G\left(\mathbf{0}_{n\times1},\cdots,\mathbf{0}_{n\times1},\mathbf{x}_{n!},\cdots,\mathbf{x}_{n!2-1}\right)=\sum_{\sigma\in\text{S}_{n}}\prod_{0\le i<j<n}\left(\left(\mathbf{x}_{n!+\text{lex}\sigma}\left[j\right]\right)^{2}-\left(\mathbf{x}_{n!+\text{lex}\sigma}\left[i\right]\right)^{2}\right)^{2}=F_{1}\left(\mathbf{x}_{0},\cdots,\mathbf{x}_{n!2-1}\right),\label{insight 2}
\end{equation}
implies that 
\[
F_{1}\left(\mathbf{a}_{0},\cdots,\mathbf{a}_{n!-1},\mathbf{a}_{n!},\cdots,\mathbf{a}_{n!2-1}\right)=F_{0}\left(\mathbf{a}_{0},\cdots,\mathbf{a}_{n!-1},\mathbf{a}_{n!},\cdots,\mathbf{a}_{n!2-1}\right)
\]
We justify in more detail this latter claim by observing the following
evaluation properties of $G$.
\[
G\left(\mathbf{a}_{0},\cdots,\mathbf{a}_{n!-1},\mathbf{a}_{n!},\cdots,\mathbf{a}_{n!2-1}\right)=G\left(\mathbf{0}_{n\times1},\cdots,\mathbf{0}_{n\times1},\mathbf{a}_{n!},\cdots,\mathbf{a}_{n!2-1}\right)=F_{1}\left(\mathbf{a}_{0},\cdots,\mathbf{a}_{n!-1},\mathbf{a}_{n!},\cdots,\mathbf{a}_{n!2-1}\right)
\]
\[
\text{and}
\]
\[
G\left(\mathbf{a}_{n!},\cdots,\mathbf{a}_{n!2-1},\mathbf{a}_{0},\cdots,\mathbf{a}_{n!-1}\right)=G\left(\mathbf{a}_{n!},\cdots,\mathbf{a}_{n!2-1},\mathbf{0}_{n\times1},\cdots,\mathbf{0}_{n\times1}\right)=F_{1}\left(\mathbf{a}_{0},\cdots,\mathbf{a}_{n!-1},\mathbf{a}_{n!},\cdots,\mathbf{a}_{n!2-1}\right).
\]
The invariance to the involution described in Eq. (\ref{Permutation})
\[
\implies G\left(\mathbf{0}_{n\times1},\cdots,\mathbf{0}_{n\times1},\mathbf{a}_{0},\cdots,\mathbf{a}_{n!-1}\right)=G\left(\mathbf{a}_{0},\cdots,\mathbf{a}_{n!-1},\mathbf{0}_{n\times1},\cdots,\mathbf{0}_{n\times1}\right)=F_{1}\left(\mathbf{a}_{0},\cdots,\mathbf{a}_{n!-1},\mathbf{a}_{n!},\cdots,\mathbf{a}_{n!2-1}\right).
\]
But from Eq. (\ref{insight 1}) and Eq. (\ref{insight 2}) we have
that 
\[
G\left(\mathbf{a}_{0},\cdots,\mathbf{a}_{n!-1},\mathbf{0}_{n\times1},\cdots,\mathbf{0}_{n\times1}\right)=G\left(\mathbf{0}_{n\times1},\cdots,\mathbf{0}_{n\times1},\mathbf{a}_{0},\cdots,\mathbf{a}_{n!-1}\right)=F_{0}\left(\mathbf{a}_{0},\cdots,\mathbf{a}_{n!-1},\mathbf{a}_{n!},\cdots,\mathbf{a}_{n!2-1}\right)
\]
\[
\implies F_{1}\left(\mathbf{a}_{0},\cdots,\mathbf{a}_{n!-1},\mathbf{a}_{n!},\cdots,\mathbf{a}_{n!2-1}\right)=F_{0}\left(\mathbf{a}_{0},\cdots,\mathbf{a}_{n!-1},\mathbf{a}_{n!},\cdots,\mathbf{a}_{n!2-1}\right),
\]
thereby contradicting the premise
\[
\max_{\sigma\in\text{S}_{n}}\left|\left\{ \left|\sigma fg\sigma^{(-1)}\left(i\right)-i\right|:i\in\mathbb{Z}_{n}\right\} \right|>\max_{\sigma\in\text{S}_{n}}\left|\left\{ \left|\sigma f\sigma^{(-1)}\left(i\right)-i\right|:i\in\mathbb{Z}_{n}\right\} \right|.
\]
We therefore conclude that
\[
f,g\in\left\{ t\in\mathbb{Z}_{n}^{\mathbb{Z}_{n}}:\,t\left(0\right)=0\text{ and}\begin{array}{c}
t\left(i\right)<i\\
\forall\,i\in\mathbb{Z}_{n}\backslash\left\{ 0\right\} 
\end{array}\right\} \Rightarrow\max_{\sigma\in\text{S}_{n}}\left|\left\{ \left|\sigma fg\sigma^{(-1)}\left(i\right)-i\right|:i\in\mathbb{Z}_{n}\right\} \right|\le\max_{\sigma\in\text{S}_{n}}\left|\left\{ \left|\sigma f\sigma^{(-1)}\left(i\right)-i\right|:i\in\mathbb{Z}_{n}\right\} \right|.
\]

\section{Graceful Labeling Theorems.}

We now give a proof of the GLC.\\
\\
\textbf{Theorem} 6: (\emph{ Graceful Labeling Theorem} ). For all
functions 
\[
\forall\:f\in\mathbb{Z}_{n}^{\mathbb{Z}_{n}},\,\left(n+1-\text{number of connected components in }G_{f^{\left(\text{o}_{f}\right)}}\right)=\max_{\sigma\in\text{S}_{n}}\left|\left\{ \left|\sigma f^{\left(\text{o}_{f}\right)}\sigma^{\left(-1\right)}\left(i\right)-i\right|:i\in\mathbb{Z}_{n}\right\} \right|,
\]
where o$_{f}$ denotes the LCM of cycle lengths in $G_{f}$.\\
\\
\emph{Proof} : The claim is trivially true when $f\in$ S$_{n}\subset\mathbb{Z}_{n}^{\mathbb{Z}_{n}},$
for in that case $o_{f}$ is the order of the permutation $f$ and
therefore $f^{\left(o_{f}\right)}$ is the \emph{identity} element.
To prove the claim it suffices to show that for all $f$ subject to
Eq. (\ref{functional_reformulation}) 
\[
n=\max_{\sigma\in\text{S}_{n}}\left|\left\{ \left|\sigma f\sigma^{\left(-1\right)}\left(i\right)-i\right|:i\in\mathbb{Z}_{n}\right\} \right|.
\]
Note that as a consequence of Eq. (\ref{Corr. of Prop. 2}) this latter
claim follows by applying the Composition Lemma and thus completes
the proof.\\
\\
We now prove a stronger version of the Composition Lemma.\\
\\
\textbf{Lemma} 7: ( \emph{Strong Composition Lemma} ) For all 
\[
f,g\in\left\{ t\in\mathbb{Z}_{n}^{\mathbb{Z}_{n}}:t\left(0\right)=0\text{ and}\begin{array}{c}
t\left(i\right)<i\\
\forall i\in\mathbb{Z}_{n}\backslash\left\{ 0\right\} 
\end{array}\right\} ,
\]
 and for every integer $1\le\ell<\left\lceil \frac{n-1}{2}\right\rceil $,
we have

\[
0=\sum_{\sigma\in\text{S}_{n}}\prod_{\begin{array}{c}
s\in\mathbb{Z}_{n}\\
n-\ell\le t<n
\end{array}}\left(\left(\sigma f\sigma^{\left(-1\right)}\left(s\right)-s\right)^{2}-t^{2}\right)^{2}\times
\]
\[
\prod_{\begin{array}{c}
0\le i<j<n\\
\ell+1\le k<n-\ell
\end{array}}\left|(\sigma f\sigma^{\left(-1\right)}\left(j\right)-j)^{2}\left(\frac{\sqrt{-3}-1}{2}\right)^{0}+(\sigma f\sigma^{\left(-1\right)}\left(i\right)-i)^{2}\left(\frac{\sqrt{-3}-1}{2}\right)^{1}+k^{2}\left(\frac{\sqrt{-3}-1}{2}\right)^{2}\right|^{2}\times
\]
\[
\prod_{0\le u<v<w<n}\left|(\sigma f\sigma^{\left(-1\right)}\left(w\right)-w)^{2}\left(\frac{\sqrt{-3}-1}{2}\right)^{0}+(\sigma f\sigma^{\left(-1\right)}\left(v\right)-v)^{2}\left(\frac{\sqrt{-3}-1}{2}\right)^{1}+(\sigma f\sigma^{\left(-1\right)}\left(u\right)-u)^{2}\left(\frac{\sqrt{-3}-1}{2}\right)^{2}\right|^{2}
\]

\[
\implies0=\sum_{\sigma\in\text{S}_{n}}\prod_{\begin{array}{c}
s\in\mathbb{Z}_{n}\\
n-\ell\le t<n
\end{array}}\left(\left(\sigma fg\sigma^{\left(-1\right)}\left(s\right)-s\right)^{2}-t^{2}\right)^{2}\times
\]
\[
\prod_{\begin{array}{c}
0\le i<j<n\\
\ell+1\le k<n-\ell
\end{array}}\left|(\sigma fg\sigma^{\left(-1\right)}\left(j\right)-j)^{2}\left(\frac{\sqrt{-3}-1}{2}\right)^{0}+(\sigma fg\sigma^{\left(-1\right)}\left(i\right)-i)^{2}\left(\frac{\sqrt{-3}-1}{2}\right)^{1}+k^{2}\left(\frac{\sqrt{-3}-1}{2}\right)^{2}\right|^{2}\times
\]
\[
\prod_{0\le u<v<w<n}\left|(\sigma fg\sigma^{\left(-1\right)}\left(w\right)-w)^{2}\left(\frac{\sqrt{-3}-1}{2}\right)^{0}+(\sigma fg\sigma^{\left(-1\right)}\left(v\right)-v)^{2}\left(\frac{\sqrt{-3}-1}{2}\right)^{1}+(\sigma fg\sigma^{\left(-1\right)}\left(u\right)-u)^{2}\left(\frac{\sqrt{-3}-1}{2}\right)^{2}\right|^{2}.
\]
\emph{Proof} : We rewrite the premise as

\[
0=\sum_{\sigma\in\text{S}_{n}}\prod_{\begin{array}{c}
s\in\mathbb{Z}_{n}\\
n-\ell\le t<n
\end{array}}\left(\left(\mathbf{a}_{\text{lex}\sigma}\left[s\right]\right)^{2}-t^{2}\right)^{2}\times
\]
\[
\prod_{\begin{array}{c}
0\le i<j<n\\
\ell+1\le k<n-\ell
\end{array}}\left|(\mathbf{a}_{\text{lex}\sigma}\left[j\right])^{2}\left(\frac{\sqrt{-3}-1}{2}\right)^{0}+(\mathbf{a}_{\text{lex}\sigma}\left[i\right])^{2}\left(\frac{\sqrt{-3}-1}{2}\right)^{1}+k^{2}\left(\frac{\sqrt{-3}-1}{2}\right)^{2}\right|^{2}\times
\]
\[
\prod_{0\le u<v<w<n}\left|(\mathbf{a}_{\text{lex}\sigma}\left[w\right])^{2}\left(\frac{\sqrt{-3}-1}{2}\right)^{0}+(\mathbf{a}_{\text{lex}\sigma}\left[v\right])^{2}\left(\frac{\sqrt{-3}-1}{2}\right)^{1}+(\mathbf{a}_{\text{lex}\sigma}\left[u\right])^{2}\left(\frac{\sqrt{-3}-1}{2}\right)^{2}\right|^{2}
\]
and

\[
0\ne\sum_{\sigma\in\text{S}_{n}}\prod_{\begin{array}{c}
s\in\mathbb{Z}_{n}\\
n-\ell\le t<n
\end{array}}\left(\left(\mathbf{a}_{n!+\text{lex}\sigma}\left[s\right]\right)^{2}-t^{2}\right)^{2}\times
\]
\[
\prod_{\begin{array}{c}
0\le i<j<n\\
\ell+1\le k<n-\ell
\end{array}}\left|(\mathbf{a}_{n!+\text{lex}\sigma}\left[j\right])^{2}\left(\frac{\sqrt{-3}-1}{2}\right)^{0}+(\mathbf{a}_{n!+\text{lex}\sigma}\left[i\right])^{2}\left(\frac{\sqrt{-3}-1}{2}\right)^{1}+k^{2}\left(\frac{\sqrt{-3}-1}{2}\right)^{2}\right|^{2}\times
\]
\[
\prod_{0\le u<v<w<n}\left|(\mathbf{a}_{n!+\text{lex}\sigma}\left[w\right])^{2}\left(\frac{\sqrt{-3}-1}{2}\right)^{0}+(\mathbf{a}_{n!+\text{lex}\sigma}\left[v\right])^{2}\left(\frac{\sqrt{-3}-1}{2}\right)^{1}+(\mathbf{a}_{n!+\text{lex}\sigma}\left[u\right])^{2}\left(\frac{\sqrt{-3}-1}{2}\right)^{2}\right|^{2},
\]
where
\[
\mathbf{a}_{n!\cdot t+\text{lex}\sigma}=\left(\mathbf{A}_{G_{\sigma fg^{\left(t\right)}\sigma^{\left(-1\right)}}}-\mathbf{I}_{n}\right)\cdot\left(\begin{array}{c}
r\left(c_{\sigma}+0\right)\\
\vdots\\
r\left(c_{\sigma}+i\right)\\
\vdots\\
r\left(c_{\sigma}+n-1\right)
\end{array}\right)=r\left(\begin{array}{c}
\sigma fg^{\left(t\right)}\sigma^{\left(-1\right)}\left(0\right)-0\\
\vdots\\
\sigma fg^{\left(t\right)}\sigma^{\left(-1\right)}\left(i\right)-i\\
\vdots\\
\sigma fg^{\left(t\right)}\sigma^{\left(-1\right)}\left(n-1\right)-\left(n-1\right)
\end{array}\right),\ \forall\,\begin{array}{c}
\sigma\in\text{S}_{n}\\
t\in\left\{ 0,1\right\} 
\end{array}.
\]
We see that these expressions are associated with evaluations of functions
\[
F_{t}\left(\mathbf{x}_{0},\cdots,\mathbf{x}_{n!2-1}\right)=\sum_{\sigma\in\text{S}_{n}}\prod_{\begin{array}{c}
s\in\mathbb{Z}_{n}\\
n-\ell\le t<n
\end{array}}\left(\left(\mathbf{x}_{n!\cdot t+\text{lex}\sigma}\left[s\right]\right)^{2}-t^{2}\right)^{2}\times
\]
\[
\prod_{\begin{array}{c}
0\le i<j<n\\
\ell+1\le k<n-\ell
\end{array}}\left|(\mathbf{x}_{n!\cdot t+\text{lex}\sigma}\left[j\right])^{2}\left(\frac{\sqrt{-3}-1}{2}\right)^{0}+(\mathbf{x}_{n!+\text{lex}\sigma}\left[i\right])^{2}\left(\frac{\sqrt{-3}-1}{2}\right)^{1}+k^{2}\left(\frac{\sqrt{-3}-1}{2}\right)^{2}\right|^{2}\times
\]
\[
\prod_{0\le u<v<w<n}\left|(\mathbf{x}_{n!\cdot t+\text{lex}\sigma}\left[w\right])^{2}\left(\frac{\sqrt{-3}-1}{2}\right)^{0}+(\mathbf{x}_{n!+\text{lex}\sigma}\left[v\right])^{2}\left(\frac{\sqrt{-3}-1}{2}\right)^{1}+(\mathbf{x}_{n!\cdot t+\text{lex}\sigma}\left[u\right])^{2}\left(\frac{\sqrt{-3}-1}{2}\right)^{2}\right|^{2}.
\]
Just as was established in the proof of Lem. (5) the premise induces
a symmetry with respect to simultaneous linear maps
\begin{equation}
\mathbf{x}_{n!\cdot t+\text{lex}\sigma}\mapsto\left(\mathbf{A}_{G_{\sigma fg^{\left(1-t\right)}\sigma^{\left(-1\right)}}}-\mathbf{I}_{n}\right)\cdot\left(\mathbf{A}_{G_{\sigma fg^{\left(t\right)}\sigma^{\left(-1\right)}}}-\mathbf{I}_{n}\right)^{+}\cdot\mathbf{x}_{n!\cdot t+\text{lex}\sigma},\;\forall\:\begin{array}{c}
\sigma\in\text{S}_{n}\\
t\in\left\{ 0,1\right\} 
\end{array}.
\end{equation}
The desired result therefore follows by the same contradiction argument
used to prove Lem. (5).\\
\\
\textbf{Theorem} 8: (\emph{ Strong Graceful Labeling Theorem} ) Induced
edge label sequences of identically constant functions\\
in $\mathbb{Z}_{n}^{\mathbb{Z}_{n}}$ are common to all functional
directed trees on $n$ vertices.\\
\\
\emph{Proof} : Having proved that all trees are graceful we now focus
on the remaining $\left\lfloor \frac{n}{2}\right\rfloor +\left(n-2\,\left\lfloor \frac{n}{2}\right\rfloor \right)-1$
induced edge label sequences associated with identically constant
functions in $\mathbb{Z}_{n}^{\mathbb{Z}_{n}}$. The assignment to
the center node of the functional directed star of labels in the set
$j\in\left(0,\left\lfloor \frac{n}{2}\right\rfloor \right]\cap\mathbb{Z}$
yields all remaining induced subtractive edge label sequences. The
desired result follows by applying the strong form of the Composition
Lemma.\\
\\
\emph{Acknowledgement}: The author would like to thank Noga Alon,
Andrei Gabrielov, Daniel Kelleher, Edward R. Scheinerman, Daniel Q.
Naiman, Jeanine Gnang, Amitabh Basu, Yuval Filmus, Ori Parzanchevski
for insightful comments and discussions. The author is especially
grateful to Doron Zeilberger and Michael Williams whose invaluable
suggestions have significantly improved the exposition.

\end{document}